\documentclass[12pt]{amsart}
\usepackage{amssymb, latexsym}
\usepackage{verbatim}
\usepackage{graphicx}
\theoremstyle{plain}
\newtheorem{theorem}{Theorem}

\newtheorem{lemma}[theorem]{Lemma}
\theoremstyle{definition}
\newtheorem{definition}[theorem]{Definition}
\newtheorem{remark}[theorem]{Remark}

\numberwithin{equation}{section}
\newcommand{\AUTHOR}[2]{\author{#1 #2}}

\begin{document}
%
%
%
%
%
\title{Dependence spaces II}
\AUTHOR {Ewa} {Graczy\'{n}ska}

\email {egracz@eranet.pl}



\maketitle

\begin{abstract}
This is a continuation of my paper \cite{10} and the presentation at
the Conference on Universal Algebra and Lattice Theory, Szeged,
Hungary, June 21-25, 2012. The Steinitz exchange lemma is a basic
theorem in linear algebra used, for example, to show that any two
bases for a finite-dimensional vector space have the same number of
elements.  The EIS property was introduced by A. Hulanicki, E.
Marczewski, E. Mycielski in \cite{EM1}  (see \S 31 of \cite{GG}). In
this note we show that its analogue holds in a dependence space.
\end{abstract}

\section{Basic notions}
According to F. G\'{e}cseg, H. J\"{u}rgensen \cite{6} the result
which is usually referred to as the "Exchange Lemma" states that for
transitive dependence, every independent set can be extended to form
a basis. In \cite{10} we discussed some interplay between notions
discussed in \cite{1}--\cite{2} and \cite{7}--\cite{6}. Another
proof was presented there, of the result of N.J.S. Hughes \cite{1}
on Steinitz' exchange theorem for infinite bases in connection with
the notions of transitive dependence, independence and dimension as
introduced in \cite{7} and \cite{8}. In that proof we assumed
Kuratowski-Zorn's Lemma (see \cite{4}), as a requirement pointed in
\cite{6}. In this note we extend the results to EIS property known
in general algebra as Exchange of Independent Sets Property.

We use a modification of the the notation of \cite{1}, \cite{2} and
\cite{6}: \\$a,b,c,...,x,y,z,...$ (with or without suffices) to
denote the elements of a space ${\bf S}$ (or {\bf S}') and
$A,B,C,...,X,Y,Z,...,$ for subsets of ${\bf S}$ (or {\bf S}').
$\Delta$, ${\it S}$ denote a family of subsets of ${\bf S}$, $n$ is
always a positive integer.

$A \cup B$ denotes the union of sets $A$ and $B$, $A+B$ denotes the
disjoint union of $A$ and $B$, $A-B$ denotes the difference of $A$
and $B$, i.e. is the set of those elements of $A$ which are not in
$B$.

\section{Dependent and independent sets}

The following definitions are due to N.J.S. Hughes, invented in 1962
in \cite{1}:

\begin{definition}
A set ${\bf S}$ is called a {\it dependence space} if there is
defined a set $\Delta$, whose members are finite subsets of ${\bf
S}$, each containing at least 2 elements, and if the Transitivity
Axiom is satisfied.
\end{definition}

\begin{definition}
A set $A$ is called {\it directly dependent} if $A \in \Delta$.
\end{definition}

\begin{definition}
An element $x$ is called {\it dependent on A} and is denoted by $x
\sim \Sigma A$ if either $x \in A$ or if there exist distinct
elements $x_{0},x_{1},...,x_{n}$ such that
\begin{center}
(1) $\{ x_{0},x_{1},...,x_{n} \} \in \Delta$
\end{center}
where $x_{0}=x$ and $x_{1},...,x_{n} \in A$

and {\it directly dependent} on $\{ x \}$ or $\{ x_{1},...,x_{n}
\}$, respectively.
\end{definition}

\begin{definition}
A set $A$ is called {\it dependent} (with respect to $\Delta$) if
(1) is satisfied for some distinct elements $x_{0},x_{1},...,x_{n}
\in A$. Otherwise $A$ is {\it independent}.
\end{definition}

\begin{definition}
If a set $A$ is  {\it independent} and for any $x \in {\bf S}$, $x
\sim \Sigma A$, i.e. $x$ is dependent on $A$, then $A$ is called a
{\it basis of} ${\bf S}$.
\end{definition}

A similar definition of a {\it dependence} $D$ was introduced in
\cite{7}. In the paper \cite{6} authors based on the theory of
dependence in universal algebras as outlined in \cite {7}.  We
accept the well known:

\begin{definition}
The {\it span} $<X>$ of a subset $X$ of ${\bf S}$ is the set of all
elements of ${\bf S}$ which depends on $X$, i.e. $x \in <X>$ iff $x
\sim \Sigma X$.
\end{definition}

\begin{definition}
{\bf TRANSITIVITY AXIOM:}
\begin{center}
If $x \sim \Sigma A$ and for all $a \in A$, $a \sim \Sigma B$, then
$x \sim \Sigma B$.
\end{center}
\end{definition}

\section{Steinitz' exchange theorem}

In Linear Algebra,  Steinitz exchange Lemma states that:

if $a \not \in <A \cup \{ b \}>$ and $b \not \in <A>$, then $b \not
\in <A \cup \{ a \}>$.

In particular, if $A$ is independent and $a \not \in <A>$, then:

$A \cup \{ a \}$ is independent.

The following lemma is a generalization of the result of P.M. Cohn
\cite{7} (cf.  the property (E) of \cite{EM3}, p.  206, called there
an axchange of an idependent sets or Theorem 3.8 of \cite{6}, p.
426):
\begin{lemma}
In a dependence space ${\bf S}$, assume that $a \not \in <A \cup \{
b \}>$ and $b \not \in <A>$. \\Then $b \not \in <A \cup \{ a \}>$.
\end{lemma}

{\it Proof}

If $b \in <A \cup \{ a \} > - <A>$, then there exists $a_1,...,a_n
\in A$, \\such that $b \sim \{a,a_1,...,a_n \}$, i.e. $\{
a,a_1,...a_n,b \} \in \Delta$.\\ Therefore $a \in < \{ b \} \cup
A>$, a contradiction . $\Box$

It is clear, that for an independent set $A$, one gets for each $a
\in A$, that $a \not \in <A - \{ a \}>$.

\begin{remark}
The Lemma above can be reformulated in a similar way as conditions
(1) and (2) in Theorem 3.8  and conditions (1) -- (3) of Lemma 3.9
of \cite{6}.
\end{remark}

\section{EIS property}
The EIS (exchange of independent sets) property was introduced by A.
Hulanicki, E. Marczewski, E. Mycielski in \cite{EM1}. First we
recall their original definition of EIS property (see \cite{EM1},
\cite{EM3}, p. 647--659). In their paper they  use the terminology
and notation of \cite{EM2} (with slight modifications).  An {\it
abstract algebra} is a (nonempty) set with a family of fundamental
finitary operations.  For  any nonempty set $E \subset A$. ${\it
C}(E)$ denotes the subalgebra  generated by $E$, ${\it
C}(\emptyset)$ is denoting the set of algebraic constants (i.e. the
values of the constant algebraic operations). The operation ${\it
C}$ has finite character, i.e. ${\it C}(E) = \bigcup {\it C}(F)$,
where $F$ runs over the family of all finite subsets of $F$ of $E$.

The following theorem about exchange of independent sets is true for
all algebras (see \cite{EM2}, p. 58, theorem 2.4 (ii)):

\begin{theorem}
Let $P,Q$ and $R$ be subsets of an algebra. If $P \cup Q$ is
independent, $P \cap Q = \emptyset$, $R$ is independent, ${\it C}(R)
= {\it C}(Q)$, then $P \cup R$ is independent.
\end{theorem}
As the authors of \cite{EM1} noticed, it might seem at first glance
that the relation ${\it C}(R) = {\it C}(Q)$ could be replaced by a
weaker one: $R \subset {\it C}(Q)$. Since, as it can be seen from
the results of \cite{EM1}, this is not generally true, the authors
say that an algebra satisfies {\it the condition of exchange of
independent sets} (EIS) whenever for any subsets $P,Q$ and $R$ of
it, the relations: $P \cup Q$ is independent, $P \cap Q =
\emptyset$, $R$ is independent and $R \subset {\it C}(Q)$ imply that
$P \cup R$ is independent.

We transform the original definition of EIS property from {\it
algebras} to {\it dependence spaces} in the natural way:
\begin{definition}
A dependence space ${\bf S}$ satisfies the EIS property, if for
arbitrary subsets $P,Q$ and $R$ of ${\bf S}$  the conditions:

(7) $P \cap Q = \emptyset$;

(8) $P \cup Q$ is an independent set in ${\bf S}$;

(9) $R$ is an independent set in ${\bf S}$, $R \subseteq <Q>$;

altogether imply that:

(10) $P \cup R$ is an independent set.

\end{definition}
\begin{theorem}
In a dependence space ${\bf S}$, the EIS property holds.
\end{theorem}

{\it Proof}

Assume (7) -- (9). \\To show (10) assume a contrario that $P \cup R$
is a dependent set. Therefore there exist (all different) elements
$a_1,...,a_n,b_1,...,b_m \in P \cup R$ with $a_1,...,a_n \in P$ and
$b_1,...,b_m \in R$ and such that $\{ a_1,...,a_n,b_1,...,b_m \}\in
\Delta$. From (7) and (9) it follows that there exists an element
$a_1 \in P$ such that $a_1 \sim \Sigma \{a_2,...,a_n,b_1,...,b_m
\}$, i.e. $a_1 \sim \Sigma ((P - \{ a_1 \}) \cup R)$. But for very
element $b \in R$, $b \sim \Sigma Q$, therefore  $b \sim \Sigma ((P
- \{ a _1 \}) \cup Q)$. Moreover, $c \sim \Sigma((P \cup Q) -
\{a_{1}\})$, for every $c \in ((P - \{a_{1}\}) \cup R)$. Thus, by
the transitivity axiom $a_1 \sim \Sigma ((P \cup Q) - \{a_1 \}))$.
That contradicts (8), as $a_1 \in P \cup Q$ and it is clear, that
for an independent set $A$, one gets for each $a \in A$, that $a
\not \in <A - \{ a \}>$. $\Box$


\end{document}